# CONSISTENCY WITHOUT CUT ELIMINATION
Kai Brünnler and Alessio Guglielmi (TU Dresden)
22.4.2002 - updated on 24.10.2002

In this note we will show how to get consistency for first order classical logic, in a purely syntactic way, without going through cut elimination. The procedure is very simple and it uses the calculus of structures [WS] in an essential way. It also shows how finitaryness (in the sense of finite choice of premises for each rule) is actually a *triviality* (contrarily to what one would guess from textbooks).

We all know that Gödel's incompleteness forbids easy solutions for very expressive cases, so this method can't be a magic bullet for the foundations of mathematics. Nonetheless, there are two points that we believe are worthy:

1) The method is general and could be helpful for easy checking of consistency in other cases.

2) This makes clear that *there is no deep connection* between cut elimination and finitaryness, as entailed by the subformula property (which in turn entails consistency).

First we give a general summary of the steps involved, for a generic logic; then we get into the details for classical logic.

## Consistency in Four Easy Steps

Step 1

Write your logic in the calculus of structures, making sure that the rules are defined according to our recipe [AG1]. This ensures that the cut rule can trivially be replaced by a cut rule in atomic form. Then, you have the atomic cut rule, which is

$$ai\uparrow \frac{S(a, \neg a)}{S\{\bot\}} \quad .$$

where $\bot$ is the unit for disjunction.

(In the sequent calculus, the cut rule cannot trivially be restricted to atomic form. You'd have to go through full-blown cut elimination.)

Step 2

Replace $ai\uparrow$ by $fai\uparrow$ (finitary cut, see [AG2]), which is:

```
        S(a,¬a)
 fai↑-------     where a or ¬a appears in S{ }.
        S{⊥}
```

Provability is not affected by this substitution: proving this claim is almost trivial, and it probably always is for every system (look at the proof below for classical logic).

This means that you get a system that proves the same logic, but all its rules can only be applied in a finite number of different ways: the equivalent of a subformula property.

This is a weak form of cut elimination: you only eliminate *certain* cuts, those that introduce atoms that have nothing to do with the conclusion of the proof. Doing this *does not require* going through each and every mutual relation of rules in the system.

Step 3

Prove that falsehood, or empty, or any given structure, can not be proved. This is consistency in a weak form, if you like. This is very easy for classical logic and BV [AG3]. We guess it's the same for other systems.

Step 4

Prove a stronger form of consistency, like: you can't have a proof of R and another of ¬R at the same time. There is a nice easy trick for doing this in the calculus of structures, given `weak' consistency.

* * *

At this point you're done. Of course, you still can perform the equivalent of the cut elimination in the sequent calculus. Why would you want to do that anyway, since you already have finitaryness and consistency? Because doing this reduces drastically the nondeterminism of the system, and because it entails several other proof theoretical properties. But it's not about consistency and finitaryness, that connection was an artifact of the sequent calculus.

## Consistency for Propositional Classical Logic

We apply the procedure above to classical propositional logic. We refer here to systems in [BT], reproduced below for convenience.

Step 1

Classical logic in the calculus of structures: here is a presentation with atomic interaction rules and non-atomic contraction and weakening. Things work the same for the entirely atomic presentation.

1) Language:

   1) R, T and U stand for generic structures;
   2) a is an atom, which is a structure;
   3) **f** and **t** are structures, but not atoms; **f** stands for the unit `false', **t** stands for the unit `true';
   4) $[R_1,…,R_h]$ and $(R_1,…,R_h)$ are structures, for n > 0; [R,T] stands for R ∨ T, (R,T) stands for R ∧ T;
   5) the structure ¬R is the negation of R;
   6) S{ } is a structure with a hole that does not occur inside a negation;
   7) S[R,T] and S(R,T) are shortcuts of S{[R,T]} and S{(R,T)}, respectively;

2) Equations:

   Structures are considered equivalent modulo =, so defined:
   1) associativity and commutativity of […] and (…): e.g., [((a,b),c),d] = [d,(a,c,b)];
   2) De Morgan: e.g., ¬[a,(b,c)] = (¬a,[¬b,¬c]).
   3) units: [**t**,**t**] = **t**;    [**f**,R] = R;
               (**f**,**f**) = **f**;   (**t**,R) = R.

3) Rules:

   1) Interaction:
   $$\mathsf{ai}{\downarrow} \frac{S\{\mathbf{t}\}}{S[a,\neg a]} \qquad \mathsf{ai}{\uparrow} \frac{S(a,\neg a)}{S\{\mathbf{f}\}}$$

   2) Core:
   $$\mathsf{s} \frac{S([R,U],T)}{S[(R,T),U]}$$

   3) Non-core:
   $$\mathsf{c}{\downarrow} \frac{S[R,R]}{S\{R\}} \qquad \mathsf{c}{\uparrow} \frac{S\{R\}}{S(R,R)}$$

   $$\mathsf{w}{\downarrow} \frac{S\{\mathbf{f}\}}{S\{R\}} \qquad \mathsf{w}{\uparrow} \frac{S\{R\}}{S\{\mathbf{t}\}}$$

A proof of R is a chain of rule instances whose premise is **t** and conclusion is R, denoted by

   $\begin{array}{c} \mathbf{t} \\ \| \\ R \end{array}$ .

## Step 2

The rules ai↑ and w↑ are non-finitary, since, given a conclusion, they yield infinitely many premises. Being non-core, w↑ can immediately be eliminated by using a generic cut and w↓:

```
     S(R,[t,f]) = S{R}
   s-----------
     S[t,(R,f)]
  w↓-----------
     S[t,(R,¬R)]
  i↑----------- .
        S{t}
```

Generic cuts can be reduced to atomic cuts by recursively applying

```
     S(R,T,[-R,-T])
   s----------------
     S(R,[-R,(T,-T)])
  i↑---------------
         S(R,-R)
      i↑------- .
          S{f}
```

So, the only infinitary rule we are left with is ai↑. How do we get rid of it in case we have a proof?

Consider the rule

```
          S(a,¬a)
  fai↑-------      where a or ¬a appears in S{ }.
          S{⊥}
```

This rule is finitary, and we're going to show that, given a proof of R, we can always get a proof of R where no ai↑ appears, but the only cuts that appear are fai↑ instances.

Take a proof of R and individuate the bottommost instance of ai↑ that violates the proviso, as in:

```
          t
          |
       S(a,¬a)
   ai↑------- ,
        S{f}
          |
          R
```

where neither a nor ¬a appear in S{ }. We can then replace all occurrences of a and ¬a in the subproof above S{f} with **t** and **f**, respectively, and we still have a proof of R. All rule instances stay valid or become trivial, for example

```
       S{t}                S(t,f)
  ai↓ ------    and   ai↑ ------  ,
       S[t,f]              S{f}
```

can just be removed, since **t** = [**t**,**f**] and (**t**,**f**) = **f**.

Please notice that if a or ¬a appeared in S{ }, this would not work, because it could destroy the derivation from S{**f**} to R.

Proceeding inductively upwards, we remove all infinitary atomic cuts.

Step 3

We have to show that there is no proof

```
    t
    | .
    f
```

After having done Step 2, we know that we can restrict ourselves to the finitary case, which means that we just have to show that there is no such proof when no atom appears in the proof. So, only f and t can appear in the proof, and we essentially have to look at the units' equations.

Let us call `boolean structure' any structure freely made from **f** and **t** by [_,_] and (_,_). It's easy to show that **f** is not equal to **t**. Then we just have to show that no rule, when applied to a boolean structure, has premise **t** and conclusion **f**. This is simply done by inspection of all rules, which means inspection of s, c↓ and w↓. Easy case analysis.

Step 4

We want to prove that if R is provable, then ¬R is not provable. The following flipping construction is not possible in the sequent calculus.

We assume that we have both proofs:

```
    t          t
    |   and    | .
    R          ¬R
```

We compose the two proofs by conjoining them:

```
    (t,t) = t
      |           .
    (R,¬R)
```

We *flip* the proof we've got, by negating everything and using the corules for every rule instance (for this we need the entire system above). We get the valid derivation

```
[¬R,R]
   |  .
   f
```

From this derivation, we can build

```
        t
 i↓------
    [¬R,R]
       |  .
       f
```

We transform the identity i↓ into several atomic identities (by using the switch inductively) and then we have a proof of **f** from **t**, which is forbidden. Contradiction.

## Consistency for Predicate Classical Logic

We only sketch the argument, by noting the differences with respect to the propositional case; see [BT] for details of the predicative system.

We eliminate all those instances of atomic cut,

$$\text{ai}\!\uparrow \frac{S(p(t_1,\ldots,t_h),\neg p(t_1,\ldots,t_h))}{S\{\mathbf{f}\}} \quad,$$

in which predicate symbol p does not occur in S.

This is done in the same way as in the propositional case. We replace all occurrences of instances of $p(t_1,\ldots,t_h)$ by **t** and of $\neg p(t_1,\ldots,t_h)$ by **f**. This works because the only extra rules in the predicative case that care about atoms are:

1) the equation

   $\exists x.R = \forall x.R = R$   if x is not free in R;

2) the rule

$$\text{n}\!\downarrow \frac{S\{R[x\leftarrow t]\}}{S\{\exists x.R\}} \quad.$$

Instances of both of them are still valid if **t** or **f** replaces arbitrary atoms in R.